\documentclass[12pt]{amsart}
 \newtheorem{thm}{Theorem}[section]
 \newtheorem{cor}[thm]{Corollary}
 \newtheorem{lem}[thm]{Lemma}
 \newtheorem{prop}[thm]{Proposition}
 \theoremstyle{definition}
 \newtheorem{defn}[thm]{Definition}
 \theoremstyle{remark}
 
 \numberwithin{equation}{section}
 \newcommand{\norm}[1]{\left\Vert#1\right\Vert}
 
 \newcommand{\C}{\mathbb{C}}

\begin{document}

\title[]
 {A Generalization of Beurling's Theorem and Quasi-Inner Functions}

\author{ Yun-Su Kim }

\address{Department of Mathematics, The University of Toledo,
2801 W. Bancroft St. Toledo, OH 43606, U.S.A.}

\email{Yun-Su.Kim@utoledo.edu}


\dedicatory{}

%\commby{Daniel J. Rudolph}

%%% ----------------------------------------------------------------------

\begin{abstract}
We introduce two kinds of quasi-inner functions. Since every
rationally invariant subspace for a shift operator $S_K$ on a
vector-valued Hardy space $H^{2}(\Omega,K)$ is generated by a
quasi-inner function, we also provide relationships of quasi-inner
functions by comparing rationally invariant subspaces generated by
them. Furthermore, we discuss fundamental properties of
quasi-inner functions, and quasi-inner divisors.
\end{abstract}

%%% ----------------------------------------------------------------------
\maketitle
%%% ----------------------------------------------------------------------

\section{Introduction}
Beurling characterized all invariant subspaces for the shift
operator on the Hardy space $H^2$ in terms of inner functions
\cite{1}. If $\varphi$ and $\phi$ are inner functions such that
$\varphi{H^{2}}\subset{\phi{H^{2}}}$, then we have $\varphi$ is
divisible by $\phi$ \cite{H}. In fact, the converse is also true
\cite{B1}.

In this paper, $\Omega$ denotes a bounded finitely connected
region in the complex plane and $R(\Omega)$ denotes the algebra of
rational functions with poles off $\overline{\Omega}$.

For a Hilbert space $K$ and a shift operator $S_{K}$ on a
vector-valued Hardy space $H^{2}(\Omega,K)$, Y.S. Kim
characterized every $R(\Omega)$-invariant subspace $M$ for the
operator $S_{K}$, that is, invariant under $u(S_{K})$ for every
$u\in{R(\Omega)}$, in terms of quasi-inner functions \cite{K};
$M=\psi{H}^{2}(\Omega,{K}^{\prime})$ for some quasi-inner function
$\psi:\Omega\rightarrow{L(K^{\prime},K)}$ and a Hilbert space
$K^\prime$. Even though a quasi-inner function is defined as an
operator-valued function in \cite{K}, by the Riesz representation
theorem, we also provide a definition of a scalar-valued
quasi-inner function.

For quasi-inner functions
$\varphi\in{H^{\infty}(\Omega,L(\C^{n}))}$ and
$u\in{H^{\infty}(\Omega)}$, we provide relationships between
operator-valued and scalar-valued quasi-inner functions (Theorem
\ref{31}). In addition, by using multiplication operator on a
vector-valued Hardy space, we characterize quasi-inner functions
(Corollary \ref{270}).

For quasi-inner functions $\theta\in{H^{\infty}(\Omega)}$ and
$\varphi\in{H^{\infty}(\Omega,K)}$, we provide definitions of the
following two cases ;

(1) $\theta$ is divisible by $\varphi$.

(2) $\varphi$ is divisible by $\theta$.

With these definitions, we characterize those divisibilities by
comparing $R(\Omega)$-invariant subspaces,
$\theta{H^{2}(\Omega,K)}$ and $\varphi{H^{2}(\Omega,K)}$.
Furthermore, the following result of Theorem \ref{34}( also
Theorem \ref{65}) is similar to that of Proposition \ref{33} which
is in case of scalar-valued inner functions on the open unit disk
; for any quasi-inner functions $\theta\in{H^{\infty}(\Omega)}$
and $\varphi\in{H^{\infty}(\Omega,L(K))}$, the following
assertions are equivalent:

(a) $\theta{|}{\varphi}$.

(b)
$\varphi{H}^{\infty}(\Omega,K)\subset\theta{H}^{\infty}(\Omega,K)$.

(c) $\varphi{H}^{2}(\Omega,K)\subset\theta{H}^{2}(\Omega,K)$.

(d) There is a $\lambda>{0}$ such that
$\varphi(z)\varphi(z)^{*}\leq{\lambda^{2}|\theta(z)|^{2}I_{K}}$
for any $z\in{\Omega}$.

The author would like to express her gratitude to her thesis
advisor, Professor Hari Bercovici.

%%% ----------------------------------------------------------------------
\section{Preliminaries and Notation}
%%% ----------------------------------------------------------------------
In this paper, $\C$, $\overline{M}$, and $L(H)$ denote the set of
complex numbers, the (norm) closure of a set $M$, and the set of
bounded linear operators from $H$ to $H$ where $H$ is a Hilbert
space, respectively.

%\section{Greatest Common Quasi-Inner Divisor}\label{o3}
\subsection{Inner functions}

Let $\textbf{D}$ be the open unit disc. We denote by $H^{\infty}$
the Banach space of all bounded analytic functions
$\phi:\textbf{D}\rightarrow{\C}$ with the norm
$\norm\phi_{\infty}=$ sup$\{\norm{\phi(z)}:z\in{\textbf{D}}\}.$

Let $\theta$ and $\theta^{\prime}$ be two functions in
$H^{\infty}$. We say that $\theta$ \emph{divides}
$\theta^{\prime}$(or $\theta$$\mid$$\theta^{\prime}$) if
$\theta^{\prime}$ can be written as
$\theta^{\prime}=\theta\cdot\phi$ for some $\phi\in {H^{\infty}}$.
We will use the notation $\theta\equiv{\theta}^{\prime}$ if
$\theta$$\mid$${\theta}^{\prime}$ and
${\theta}^{\prime}$$\mid$${\theta}$.

 Recall that a function
$u\in{H^\infty}$ is \emph{inner} if $\mid$$u(e^{it})$$\mid$$=1$
almost everywhere on $\partial{\textbf{D}}$.
 By Beurling's theorem on
invariant subspaces of the Hardy spaces, for any inner function
$\theta\in{H^{\infty}}$, we have that $\theta{H}^{2}$ is an
invariant subspace for the shift operator
$S:H^{2}\rightarrow{H^{2}}$ defined by $(Sf)(z)=zf(z)$ for
$f\in{H^2}$.
 Generally, the
following Proposition shows how to define a \emph{divisor} of an
inner function in $H^{\infty}$ by using $S$-invariant subspaces.

\begin{prop}(\cite{B1})\label{33}
For any inner functions $\theta$ and $\theta^\prime$ in
$H^{\infty}$, the following assertions are equivalent:

(a) $\theta\mid{\theta^{\prime}}$.

(b) $\theta^{\prime}H^{\infty}\subset\theta{H}^{\infty}$.

(c) $\theta^{\prime}H^{2}\subset\theta{H}^{2}$.

(d) $\mid$$\theta^{\prime}(z)$$\mid\leq{c}\mid$$\theta(z)$$\mid$
for some $c>{0}$ and for all $z\in{\textbf{D}}$.
\end{prop}

\subsection{Hardy spaces}
We refer to \cite{R2} for basic facts about Hardy space, and
recall here the basic definitions. Let $\Omega$ be a bounded
finitely connected region in the complex plane.
\begin{defn}
The space ${\: H^{2}(\Omega)}$ is defined to be the space of
analytic functions $f$ on $\Omega$ such that the subharmonic
function $|f|^{2}$ has a harmonic majorant on $\Omega$. For a
fixed $z_{0}$ $\in\Omega$, there is a norm on
$H^{2}$$($$\Omega$$)$ defined by

      $\|f\|$=inf$\{ u(z_{0})^{1/2}$: $u$ is a harmonic majorant of
      $|f|^{2}\}$.
\end{defn}

 Let $m$ be the harmonic measure for the point $z_{0}$, let
 $L^{2}(\partial{\Omega})$ be the $L^{2}$-space of complex valued
 functions on the boundary of $\Omega$ defined with respect to $m$,
 and let $H^{2}(\partial{\Omega})$ be the set of
 functions $f$ in $L^{2}(\partial{\Omega})$ such that
 $\int_{\partial{\Omega}} f(z)g(z) dz$ = 0 for every $g$ that is
 analytic in a neighborhood of the closure of $\Omega$.
 If $f$ is in $H^{2}(\Omega)$, then there is a function $f^{\ast}$ in
 $H^{2}(\partial{\Omega})$ such that $f({z})$ approaches
 $f^{\ast}(\lambda_{0})$ as $z$ approaches $\lambda_{0}$
 nontangentially, for almost every $\lambda_{0}$ relative to $m$. The map
 $f\rightarrow{f^{\ast}}$ is an isometry from $H^{2}(\Omega)$ onto
 $H^{2}(\partial{\Omega})$. In this way, $H^{2}(\Omega)$ can be
 viewed as a closed subspace of $L^{2}(\partial{\Omega})$.

 A function $f$ defined on $\Omega$ is in $H^{\infty}(\Omega)$ if
 it is holomorphic and bounded. $H^{\infty}(\Omega)$ is a closed
 subspace of $L^{\infty}({\Omega})$ and it is a Banach
 algebra if endowed with the supremum norm. Finally, the mapping
 $f\rightarrow{f^{\ast}}$ is an isometry of $H^{\infty}(\Omega)$
 onto a week$^{*}$-closed subalgebra of
 $L^{\infty}(\partial{\Omega})$.

\begin{defn} If $K$ is a Hilbert space, then $H^{2}(\Omega$,$K)$ is
defined to be the space of analytic functions
$f:\Omega\rightarrow{K}$ such that the subharmonic function
${\norm{f}}^{2}$ is majorized by a harmonic function $\nu$. Fix a
point $z_{0}$ in $\Omega$ and define a norm on $H^{2}(\Omega$,$K)$
by

$\norm{f}$=inf $\{{\nu(z_{0})}^{1/2}$ : $\nu$ is a harmonic
majorant of ${\norm{f}}^{2}\}$.
\end{defn}

We will work on this vector-valued Hardy space
$H^{2}(\Omega$,$K)$. Note that $H^{2}(\Omega,K)$ can be identified
with a closed subspace of the space $L^{2}(\partial\Omega$,$K)$ of
square integrable $K$-valued functions on $\partial\Omega$. Define
a shift operator
$S_{K}:H^{2}(\Omega,K)\rightarrow{H^{2}(\Omega,K)}$ by
$(S_{K}f)(z)=zf(z)$.

\section{Quasi-Inner Functions}
Let $R(\Omega)$ denote the algebra of rational functions with
poles off $\overline{\Omega}$, and $T$ be an operator in $L(H)$
such that $\sigma(T)\subset\overline{\Omega}$. Then a closed
subspace $M$ is said to be $R(\Omega)$\emph{-invariant}
(\emph{rationally invariant}) for the operator $T$, if it is
invariant under $u(T)$ for any function $u\in{R(\Omega)}$.

To characterize every $R(\Omega)$-invariant subspace for the shift
operator $S_{K}$, quasi-inner function was defined in \cite{K}.

\begin{defn}\label{26}Let $K$ and $K^{\prime}$ be Hilbert spaces and let
${H}^{\infty}(\Omega,L(K,K^{\prime}))$ be the Banach space of all
analytic functions $\Phi:\Omega\rightarrow{L(K,K^{\prime})}$ with
the supremum norm. For
$\varphi\in{H}^{\infty}(\Omega,L(K,K^{\prime}))$, we will say that
$\varphi$ is \emph{quasi-inner} if there exists a constant $c>0$
such that
\[\norm{\varphi(z)k}\geq{c\norm{k}}\]
for every $k\in{K}$ and almost every $z\in\partial\Omega$.
\end{defn}

Even though a quasi-inner function is defined as an
operator-valued function, by the Riesz representation theorem, we
can identify $L(\C)$ with $\C$. Thus we have the following
definition of a scalar-valued quasi-inner function;

\begin{defn}\label{27}
For $\theta\in{H^{\infty}(\Omega)}$, we will say that $\theta$ is
\emph{quasi-inner} if there exists a constant $c>0$ such that
\[|\theta(z)|\geq{c}\]
for almost every $z\in\partial\Omega$.

\end{defn}

\begin{prop}\label{30}
Let $K$ and $K^{\prime}$ be Hilbert spaces with $\dim K= \dim
{K}^{\prime}= n(<\infty)$.

If $\varphi\in{H}^{\infty}(\Omega,L(K,K^{\prime}))$ is a
quasi-inner function, then $\varphi(z)$ is invertible a.e. on
$\partial\Omega$.
\end{prop}
\begin {proof}

Since $\varphi\in{H}^{\infty}(\Omega,L(K,K^{\prime}))$, $\varphi$
has nontangential limits for any $z$ $\in$
$\partial\Omega\setminus$$A$ with $m(A)=0$. For a fixed $z_{0}$
$\in$ $\partial\Omega\setminus$$A$, since for some $c>0$ and $a$
$\in$ $K$,
\begin{center}$c\norm{a}\leq\norm{\varphi(z_{0})a}$,\end{center} the range of
$\varphi(z_{0})$ is closed, and $\varphi(z_{0})$ is one-to-one. By
the first isomorphism theorem, \begin{center}$\varphi(z_{0})K$
$\cong$ the range of $ \varphi(z_{0})$.\end{center} Thus, $K$ and
the range of $ \varphi(z_{0})$ have the same dimension.

Since $\dim$ $K$ = $\dim$ $K^{\prime}$ and the range of $
\varphi(z_{0})$ is a closed subspace of $K^{\prime}$, we conclude
that the range of $ \varphi(z_{0})$ is $K^{\prime}$, that is,
$\varphi(z_{0})$ is one-to-one and onto. Thus $\varphi(z)$ is
invertible for $z$ $\in$ $\partial\Omega\setminus$$A$.\end{proof}

\begin{cor}\label{48}
Let $K$ and $K^{\prime}$ be Hilbert spaces with dim $K=$ dim
$K^{\prime}=$ n$(<\infty)$. If
$\varphi\in{H}^{\infty}(\Omega,L(K,K^{\prime}))$ is quasi-inner,
and $\psi\in{H}^{\infty}(\Omega,L(K^{\prime},K))$ such that
$\varphi(z)\psi(z)=
 u(z)I_{K^{\prime}}$ for $u\in{{H}^{\infty}}(\Omega)$, $u\neq{0}$,
 then $\psi(z)\varphi(z)=
 u(z)I_{K}$.
\end{cor}
\begin{proof}
By Proposition \ref{30} (a), $\varphi(z)$ is invertible a.e. on
$\partial\Omega$. From the equation $\varphi(z)\psi(z)=
 u(z)I_{K^{\prime}}$, we obtain \begin{center}$\psi(z)=
 \varphi(z)^{-1}[u(z)I_{K^{\prime}}]=u(z)\varphi(z)^{-1}I_{K^{\prime}}=u(z)I_{K}\varphi(z)^{-1}$\end{center}
for a.e.  $z\in\partial\Omega$. It follows that for a.e.
$z\in\partial\Omega$, $\psi(z)\varphi(z)=
 u(z)I_{K}$ which proves this Corollary.
\end{proof}
In the next theorem, we have another results of relationships
between operator-valued and scalar-valued quasi-inner functions.

\begin{thm}\label{31}
(a) If $\varphi\in{H^{\infty}(\Omega,L(\C^{n}))}$ and
$u\in{H^{\infty}(\Omega)}$ are quasi-inner functions such that
\begin{center}$\varphi(z)\psi(z)=u(z)I_{\C^n}$\end{center} where
$\psi\in{H^{\infty}(\Omega,L(\C^{n}))}$, then $\psi$ is also
quasi-inner. \vskip0.1cm (b) Conversely, if
$\varphi\in{H^{\infty}(\Omega,L(\C^{n}))}$ and
$\psi\in{H^{\infty}(\Omega,L(\C^{n}))}$ are quasi-inner functions
such that \begin{center}$\varphi(z)\psi(z)=u(z)I_{\C^n}$ or
$\psi(z)\varphi(z)=u(z)I_{\C^n}$\end{center} for some
$u\in{H^{\infty}(\Omega)}(u\neq{0})$, then u is quasi-inner.
\end{thm}
\begin{proof}
(a) Since $\varphi$ and $u$ are quasi-inner functions, there is
$m_{1}>{0}$ such that $m_{1}\leq{|}u(z){|}$ a.e.
$z\in{\partial\Omega}$, and there exists $c_{i}(>0)(i=1,2)$ such
that for $h\in{\C^{n}}$,
$c_{1}\norm{h}\leq\norm{\varphi(z)h}\leq{c_{2}}\norm{h}$ a.e.
$z\in\partial\Omega$. Since $\varphi(z)\psi(z)=u(z)I_{C^{n}}$, for
$h\in{\C^{n}}$, $m_{1}\norm{h}\leq$
$|{u}(z)|\norm{h}=\parallel$$\varphi(z)\psi(z)h$$\parallel\leq{c_{2}}\norm{{\psi}(z)h}$.
Thus, for $h\in{\C^{n}}$,
\begin{equation}\label{84} \frac{m_{1}}{c_{2}}\norm{h}\leq\norm {\psi(z)h}\end{equation} a.e. $z\in{\partial\Omega}$.

 %Next, ${c_{1}}\norm{
%\psi(z)h}\leq\parallel$$\varphi(z)\psi(z)h$$\parallel=|u(z)|\norm{h}\leq\norm{u}_{\infty}\norm{h}$
%a.e. $z\in{\partial\Omega}$. Thus
%\begin{equation}\label{85}
%\norm {\psi(z)h}\leq\frac{\norm{u}_{\infty}}{c_{1}}\norm{h}
%\end{equation} a.e. $z\in{\partial\Omega}$.
From (\ref{84}), we conclude that $\psi$ is also quasi-inner.

 \vskip0.1cm (b) Since
$\varphi\in{H^{\infty}(\Omega,L(\C^{n}))}$ and
$\psi\in{H^{\infty}(\Omega,L(\C^{n}))}$ are quasi-inner functions,
there exist $m_{1}(>0)$, and $m_{2}(>0)$ such that for
$h\in{\C^{n}}$, $m_{1}\norm{h}\leq\norm{\varphi(z)h}$ and
$m_{2}\norm{h}\leq\norm{\varphi(z)h}$ a.e. $z\in\partial\Omega$.
Then
\begin{center}$\norm{\varphi(z)\psi(z)h}\geq{m_1}\norm{\psi(z)h}\geq{m_{1}m_{2}\norm{h}}$\end{center}
and so $\norm{\varphi(z)\psi(z)}\geq{m_{1}m_{2}}$. Let
$m^{\prime}=m_{1}m_{2}$. Since
$\mid$$u(z)$$\mid=\norm{\varphi(z)\psi(z)}$, it is proven.
\end{proof}

\begin{prop}\label{49}
Let $\varphi\in$ $H^{\infty}(\Omega,L(K,K^{\prime}))$. If
$\varphi$ is a quasi-inner function, then
$\varphi{H}^{2}(\Omega,K)$ is closed subspace of
$H^{2}(\Omega,K^{\prime})$.
\end{prop}
\begin{proof}
Clearly $\varphi{H}^{2}(\Omega,K)\subset{H}^{2}(\Omega,
K^{\prime})$. Let $\{\varphi{f}_{n}\}$ be a Cauchy sequence in
$\varphi{H}^{2}(\Omega,K)$. If $f_{n}^{\ast}$ is the nontangential
limits of $f_{n}$, then since $\varphi$ is a quasi-inner function,
there exists a constant $c>0$ such that
\begin{equation}\label{50} \norm{\varphi(z)f_{n}^{\ast}(z) -
\varphi(z)f_{m}^{\ast}(z)}\geq c \norm{f_{n}^{\ast}(z) -
f_{m}^{\ast}(z)} \end{equation} for a.e. $z$ $\in$
$\partial\Omega$.

Since $\norm{\varphi{f}_{n}^{\ast}}_{H^{2}(\partial\Omega,
K)}=\norm{\varphi{f}_{n}}_{H^{2}(\Omega, K)}$, $\{
\varphi{f}_{n}^{\ast}\}$ is a Cauchy sequence in
$\varphi{H}^{2}(\partial\Omega,K)$, from the equation (\ref{50})
we conclude that $\{ f_{n}^{\ast}\}$ is also a Cauchy sequence in
$H^{2}(\partial\Omega, K)$. We know that $H^{2}(\partial\Omega,
K)$ is a Banach space. Thus
$\lim_{n\rightarrow\infty}f_{n}^{\ast}$ exists.

If $\lim_{n\rightarrow\infty}
f_{n}^{\ast}={f}^{\ast}\in{H}^{2}(\partial\Omega, K)$, then there
is unique $f\in{H}^{2}(\Omega,{K})$ such that
$\lim_{n\rightarrow\infty} f_{n}=f$ and so
$\lim_{n\rightarrow\infty}\varphi{f}_{n}=\varphi{f}\in
\varphi{H}^{2}(\Omega,K^{\prime})$ which proves this Proposition.
\end{proof}

Furthermore, by using these quasi-inner functions, we have a
generalization of Beurling's theorem as following:

\noindent\textbf{Theorem A}. (Theorem 1.5 in \cite{K}) Let $K$ be
a Hilbert space. Then a closed subspace $M$ of ${H}^{2}(\Omega,
K)$ is $R(\Omega)$-invariant for $S_{K}$ if and only if there is a
Hilbert space $K^{\prime}$ and a quasi-inner function
$\varphi:\Omega\rightarrow{L(K^{\prime},K)}$ such that
$M=\varphi{H}^{2}(\Omega,K^{\prime})$.

Also, we compare two $R(\Omega)$-invariant subspaces for $S_{K}$.
\begin{prop} \cite{K}\label{23} Let $\varphi_{1}
:\Omega\rightarrow{L(K_{1},K)}$ and
$\varphi_{2}:\Omega\rightarrow{L(K_{2},K)}$ be quasi-inner
functions.

Then the subspaces $\varphi_{1}H^{2}(\Omega,K_{1})$ and
$\varphi_{2}H^{2}(\Omega,K_{2})$ of $H^{2}(\Omega,K)$ are equal if
and only if there exist functions
 $\varphi\in{H}^{\infty}(\Omega,L(K_{1},
K_{2}))$ and $\psi\in{H}^{\infty}(\Omega,L(K_{2},K_{1}))$ such
that $\varphi\psi=I_{K_{2}}$, $\psi\varphi=I_{K_{1}}$ and
$\varphi_{1}(z)=\varphi_{2}(z)\varphi(z)$ for any $z\in\Omega$.
\end{prop}

Since we have two kinds of quasi-inner functions, we have two
kinds of $R(\Omega)$-invariant subspaces for $S_{K}$. One of them
is generated by a scalar-valued quasi-inner function, and the
other one is generated by an operator-valued quasi-inner function.
We will also compare these two $R(\Omega)$-invariant subspaces for
$S_{K}$ in Theorem \ref{34}.

Let $K_1$ and $K_2$ be separable Hilbert spaces. To characterize
quasi-inner functions, we define a \emph{multiplication operator}
for a given function $\psi\in{H^{\infty}}(\Omega,L(K_{1},K{2}))$.
A \emph{multiplication operator}
$M_{\psi}:H^{2}(\Omega,K_{1})\rightarrow{H^{2}(\Omega,K_{2})}$ is
defined by
\[M_{\psi}(g)(z)=\psi(z)g(z)\]
for all $g$ in $H^{2}(\Omega,K_{1})$. We can easily check that
$\norm{M_{\psi}}=\norm{\psi}_{\infty}$.

Recall an important property of this multiplication operator :
\begin{prop}\label{38} \cite{K}
Let $K_{1}$ and $K_{2}$ be separable Hilbert spaces. If
$T:{H}^{2}(\Omega,K_{1})\rightarrow{H}^{2}(\Omega,K_{2})$ is a
bounded linear operator such that $TS_{K_{1}}=S_{K_{2}}T$, then
there is a function $\psi\in$$H^{\infty}(\Omega,L(K_{1},K_{2}))$
such that $T=M_{\psi}$, where $M_{\psi}(g)(z)=\psi(z)g(z)$ for
$g\in$${H}^{2}(\Omega,K_{1})$, we have
$\norm{T}=\norm{\psi}_{\infty}$.
\end{prop}

%A bounded operator on $H^{2}(\Omega,K)$ that commutes with the
%shift operator $S_{K}$ must be a multiplication operator \cite{K}.
\begin{cor}\label{270}
Let $\varphi\in{H^{\infty}(\Omega,L(K_{1},K_{2}))}$.

(a) If $\varphi$ is quasi-inner, then $M_{\varphi}$ is one-to-one
and has closed range.

(b) If
$M_{\varphi}:H^{2}(\Omega,K_{1})\rightarrow{H^{2}(\Omega,K_{2})}$
is invertible, then $\varphi$ is quasi-inner.
\end{cor}
\begin{proof}
$(a)$ By Proposition \ref{49},
$M_{\varphi}H^{2}(\Omega,K_{1})=\varphi{H^{2}(\Omega,K_{1})}$ is
closed. Since $\varphi$ is quasi-inner, for some $c>0$,
$\norm{\varphi(z)a}\geq{c\norm{a}}$ for a.e. $z\in\partial\Omega$
and for any $a\in{K_{1}}$. Thus condition $(a)$ implies
$\norm{M_{\varphi}f}\geq{c\norm{f}}$ for
$f\in{H^{2}(\Omega,K_{1})}$. It follows that $M_{\varphi}$ is
one-to-one. \vskip0.2cm

 $(b)$ Since
$\varphi{H^{2}(\Omega,K_{1})}$ is is $R(\Omega)$-invariant for
$S_{K_2}$, by Theorem \textbf{A},
\begin{center}$\varphi{H^{2}(\Omega,K_{1})} =\varphi_{1}{H}^{2}(\Omega,K_{0})$\end{center}
for a Hilbert space $K_{0}$ and a quasi-inner function
$\varphi_{1}:\Omega\rightarrow{L(K_{0},K_{2})}$.

Define an operator
$T:{H^{2}(\Omega,K_{1})}\rightarrow{H^{2}(\Omega,K_{0})}$ as
follows. For $f\in{H^{2}(\Omega,K_{1})}$, $Tf=g$ such that
$\varphi{f}=\varphi_{1}g$. Since $\varphi_{1}$ is a quasi-inner
function, by $(a)$, $T$ is well-defined. Since
$S_{K_{0}}T=TS_{K_{1}}$, by Proposition \ref{38},
$T=M_{\varphi_2}$ for a function
$\varphi_{2}\in$$H^{\infty}(\Omega,L(K_{1},K_{2}))$. It follows
that
\begin{equation}\label{37}\varphi(z)=\varphi_{1}(z)\varphi_{2}(z)\end{equation} for any $z\in\Omega$.

Since
$M_{\varphi}:H^{2}(\Omega,K_{1})\rightarrow{H^{2}(\Omega,K_{2})}$
is invertible, so is $T=M_{\varphi_2}$. Note that the
invertibility of $T$ is equivalent to the invertibility of
$\varphi_{2}(z)$ for any $z$ in $\Omega$.  It follows that
$\varphi_{2}(z)$ is bounded below for any $z$ in $\Omega$.

%By the same way as the proof of previous Proposition in \cite{K},
%we can get an invertible function
%$\varphi_{2}\in{H^{\infty}(\Omega,L(K_{1},K_{0}))}$ such that
%$\sup \norm{\varphi_{2}(z)^{-1}}<\infty$ and
%\begin{center}$\varphi(z)=\varphi_{1}(z)\varphi_{2}(z)$\end{center} for any $z\in\Omega$.

%Let $c_{1}=\sup \norm{\varphi_{2}(z)^{-1}}>0$. Since
%$1\leq\norm{\varphi_{2}(z)}\norm{\varphi_{2}(z)^{-1}}$,
%\begin{equation}\label{01}
%\frac{1}{c_1}\leq\frac{1}{\norm{\varphi_{2}(z)^{-1}}}\leq\norm{\varphi_{2}(z)}.
%\end{equation}
%From the fact that $\varphi_{1}$ is quasi-inner, we can conclude
%that there is a $c_{2}>0$ such that for any $a\in{K_{1}}$,
%\begin{equation}\norm{\varphi(z)a}=\norm{\varphi_{1}(z)\varphi_{2}(z)a}\geq{c_{2}}\norm{\varphi_{2}(z)a}
%\geq\frac{c_2}{c_1}\norm{a},\end{equation} for
%a.e.$z\in\partial\Omega$.
Since $\varphi_1$ is quasi-inner, $\varphi_1(z)$ is also bounded
below almost every $z\in{\partial\Omega}$.

Therefore, by equation (\ref{37}), for any $a\in{K_{1}}$, there is
a constant $c>{0}$ such that
\begin{center}$\norm{\varphi(z)a}\geq{c\norm{a}}$\end{center} for
a.e.$z\in\partial\Omega$.
\end{proof}

\section{Quasi-inner Divisors}

Let $K$ be a Hilbert space. The time has come to consider
divisibilities between a function in $H^{\infty}(\Omega)$ and a
function in $H^{\infty}(\Omega,L(K))$ or between functions in
$H^{\infty}(\Omega,L(K))$.
\begin{defn}\label{55}If $\theta\in{H^{\infty}(\Omega)}$ and
$\varphi\in{H^{\infty}(\Omega,L(K))}$, then we say that $\theta$
\emph{divides} $\varphi$ (denoted $\theta|\varphi)$ if $\varphi$
can be written as
\begin{center}$\varphi=\theta\cdot\phi^{\prime}$\end{center} for
some $\phi^{\prime}\in{H^{\infty}(\Omega,L(K))}$.\end{defn}

\begin{defn}\label{86}If $\theta\in{H^{\infty}(\Omega)}$ and
$\varphi\in{H^{\infty}(\Omega,L(K))}$, then we say that $\varphi$
\emph{divides} $\theta$ (denoted $\varphi|\theta)$ if there exists
$\psi\in{H}^{\infty}(\Omega,L(K))$ satisfying the following
relations $;$
\begin{center}$\varphi(z)\psi(z)=
 \theta(z)I_{K}$\end{center}  and \begin{center}$\psi(z)\varphi(z)=
 \theta(z)I_{K}$\end{center} for $z\in\Omega$.\end{defn}

\begin{defn}\label{56}
If $\varphi$ and $\varphi^{\prime}$ are functions in
$H^{\infty}(\Omega,L(K))$, then we say that $\varphi$ is a
\emph{left divisor} of $\varphi^\prime$ if
\begin{center}$\varphi^{\prime}(z)=\varphi(z)\varphi^{\prime\prime}(z)$ (denoted
$\varphi{|}_{l}\varphi^\prime$)\end{center} for some
$\varphi^{\prime\prime}\in{H^{\infty}(\Omega,L(K))}$, and
 we say that $\varphi$
is a \emph{right divisor} of $\varphi^\prime$ if
\begin{center}$\varphi^{\prime}(z)=\varphi^{\prime\prime}(z)\varphi(z)$ (denoted
$\varphi{|}_{r}\varphi^\prime$)\end{center} for some
$\varphi^{\prime\prime}\in{H^{\infty}(\Omega,L(K))}$. We will say
$\varphi$ \emph{divides} $\varphi^\prime$ if $\varphi$ is not only
a left divisor but also a right divisor of $\varphi^\prime$.
\end{defn}

\begin{defn}\label{58}
Let $\varphi$ be a functions in $H^{\infty}(\Omega,L(K))$. A
quasi-inner function $\theta\in{H^{\infty}(\Omega)}$ is called the
\emph{greatest quasi-inner divisor} of $\varphi$ if $\theta$
divides $\varphi$ and every (complex-valued) quasi-inner divisors
of $\varphi$ is a divisor of $\theta$. \emph{The greatest
quasi-inner divisor} of $\varphi$ is denoted by $D(\varphi)$.
\end{defn}

With definitions of divisibility, we provide the following results
similar to Proposition \ref{33}.

\begin{thm}\label{34}
For any quasi-inner functions $\theta\in{H^{\infty}(\Omega)}$ and
$\varphi\in{H^{\infty}(\Omega,L(K))}$, the following assertions
are equivalent:\vskip0.2cm

(a) $\theta{|}{\varphi}$.

(b)
$\varphi{H}^{\infty}(\Omega,K)\subset\theta{H}^{\infty}(\Omega,K)$.

(c) $\varphi{H}^{2}(\Omega,K)\subset\theta{H}^{2}(\Omega,K)$.

(d) There is a $\lambda>{0}$ such that
$\varphi(z)\varphi(z)^{*}\leq{\lambda^{2}|\theta(z)|^{2}I_{K}}$
for any $z\in{\Omega}$, where $I_K$ is the identity function on
$K$.

\end{thm}
\begin{proof}
If $\theta{|}{\varphi}$, $\varphi=\theta\varphi_{1}$ for some
$\varphi_{1}\in{H^{\infty}(\Omega,L(K))}$. Then
\begin{center}$\varphi{H}^{\infty}(\Omega,K)=\theta\varphi_{1}{H}^{\infty}(\Omega,K)
\subset\theta{H}^{\infty}(\Omega,K)$.\end{center} Thus (a) implies
(b).

Conversely, suppose that
$\varphi{H}^{\infty}(\Omega,K)\subset\theta{H}^{\infty}(\Omega,K)$.
Then
\[\varphi^{\ast}{H}^{\infty}(\partial\Omega,K)\subset\theta^{\ast}{H}^{\infty}(\partial\Omega,K).\]

Let $\{b_{i}:i\in{I}\}$ be an orthonormal basis of $K$ and
$g_{i}\in{H}^{\infty}(\partial\Omega,K)$ defined by
$g_{i}(z)=b_{i}(i\in{I})$. Since
$\varphi^{\ast}{H}^{\infty}(\partial\Omega,K)\subset\theta^{\ast}{H}^{\infty}(\partial\Omega,K)$,
there is $f_{i}\in{{H}^{\infty}(\partial\Omega,K)}$ such that
$\varphi^{\ast}{g_{i}}=\theta^{\ast}{f_{i}}$, i.e. for $i\in{I}$,
\begin{equation}\label{66}\varphi^{\ast}(z)b_{i}=\theta^{\ast}(z)f_{i}(z).\end{equation}
Define $\varphi_{1}:\partial\Omega\rightarrow{L(K)}$ by for
$i\in{I}$,
\begin{equation}\label{67}\varphi_{1}(z)b_{i}=f_{i}(z).\end{equation}
For $i\in{I}$, define
$\varphi_{i}\in{H^{\infty}(\partial\Omega,L(K))}$ by
$\varphi_{i}(z)b_{j}=\delta_{ij}f_{i}(z) (j\in{I})$, where
$\delta_{ij}=\left\{\begin{array}{ll}
                                         1 & \mbox{if $i=j$}\\0 &
                                         \mbox{otherwise.}\end{array}
                                  \right. $ Then
                                  \begin{equation}\label{b2}\varphi_{1}=\sum_{i\in{I}}\varphi_{i}.\end{equation}

\noindent By (\ref{66}) and (\ref{67}), for each $i\in{I}$,
$\varphi^{\ast}(z)b_{i}=\theta^{\ast}(z)\varphi_{1}(z)b_{i}$, and
so
\begin{center}$\varphi^{\ast}=\theta^{\ast}\varphi_{1}$.\end{center}

 To prove that $(b)$
implies $(a)$, we have to show that
$\varphi_{1}\in{H^{\infty}(\partial\Omega,L(K))}$. Since
$\theta\in{H^{\infty}}$ is a quasi-inner function, there is $c>0$
such that $|\theta(z)|\geq{c}$ for every
$z\in{A}\subset\partial\Omega$ with $m(A)=0$. For any $x\in{K}$
with $\norm{x}=1$ and $z\in{A}$,
\begin{equation}\label{b3}\norm{\varphi_{1}(z)x}=\frac{\norm{\varphi^{\ast}(z)x}}{|\theta^{\ast}(z)|}\leq
\frac{\norm{\varphi}_{\infty}}{c}.\end{equation} From (\ref{b2})
and (\ref{b3}), we conclude that
\begin{center}$\varphi_{1}\in{H^{\infty}(\partial\Omega,L(K))}$.\end{center}

  By the same way as above,
$(a)\Leftrightarrow(c)$ is proven. We begin to prove
$(a)\Leftrightarrow(d)$. If $\theta{|}{\varphi}$,
$\varphi=\theta\varphi_{1}$ for some
 $\varphi_{1}\in{H^{\infty}(\Omega,L(K))}$. Then
\begin{center}$\varphi(z)\varphi(z)^{*}=\theta(z)\varphi_{1}(z)$
$\varphi_{1}(z)^{*}\overline{\theta(z)}
\leq{\norm{\varphi_{1}}_{\infty}^{2}|\theta(z)|^{2}I_{K}}$.\end{center}
Let $\lambda=\norm{\varphi_{1}}_{\infty}$. Since $\varphi$ is
quasi-inner, $\varphi\neq{0}$ and so $\lambda>0$. Thus (a) implies
(d).\vskip0.2cm

Conversely, suppose that for any $z\in{\Omega}$,
\begin{equation}\label{69}\varphi(z)\varphi(z)^{*}\leq{\lambda^{2}|\theta(z)|^{2}I_{K}}\end{equation}
for some $\lambda>{0}$. For each $z\in{\Omega}$, we will define a
linear mapping $F_{z}\in{L(K)}$. Let
\begin{center}$A=\{z\in{\Omega}:\theta(z)=0\}$\end{center} and
\begin{center}$B=\{z\in{\Omega}:\theta(z)\neq0\}$.\end{center} If $z\in{A}$, then let
$F_{z}=0$. If $z\in{B}$, then range of $\overline{\theta(z)}I_{K}$
is $K$ and so we can define a linear mapping $F_{z}$ from $K$ to
range
 of $\varphi(z)^{*}$ by \begin{center}$F_{z}(\overline{\theta(z)}f)=\varphi(z)^{*}f$\end{center}
for $f\in{K}$.

Since
$\norm{F_{z}(\overline{\theta(z)}f)}^{2}=\norm{\varphi(z)^{*}f}^{2}
=(\varphi(z)\varphi(z)^{*}f,f)\leq\lambda^{2}(|\theta(z)|^{2}f,f)=\lambda^{2}\norm{\theta(z)f}^{2}$,
that is,
\begin{equation}\label{70}\norm{F_{z}(\overline{\theta(z)}f)}\leq{\lambda\norm{\theta(z)f}},\end{equation}
$F_{z}$ is well-defined for $z\in{B}$. By definition of $F_{z}$,
if $z\in{B}$,
\begin{equation}\label{q}\theta(z)F_{z}^{\ast}=\varphi(z).\end{equation} If
$z\in{A}$, by (\ref{69}) $\norm{\varphi(z)}=0$ and so
$\varphi(z)=0(z\in{A})$. Thus $\theta(z)F_{z}^{*}=\varphi(z)$ for
any $z\in{\Omega}$.

Define a function $F:\Omega\rightarrow{K}$ by
\begin{center}$F(z)=F_{z}^{*}$.\end{center} Then by equation (\ref{q}), \begin{center}$\varphi(z)=\theta(z)F(z)$\end{center}
for $z\in{\Omega}$. To finish this proof, we have to prove that
$F\in{H^{\infty}(\Omega,L(K))}$. From inequality (\ref{70}), we
have
\begin{equation}\label{c3}\norm{F}_{\infty}\leq{\lambda}\end{equation} and so
$F=\frac{\varphi}{\theta}$ has only removable singularities in
$\Omega$. Thus $F$ can be defined on
$\{z\in{\Omega:\theta(z)=0}\}$ so that $F$ is analytic and
\begin{center}$\varphi=\theta{F}$.\end{center} From (\ref{c3}),
$F\in{H^{\infty}(\Omega,L(K))}$ which proves $(d)\Rightarrow(a)$.
\end{proof}

Note that for $\varphi\in{H^{\infty}(\Omega,L(\C^{n}))}$, we can
think $\varphi(z)(z\in{\Omega})$ as an $n\times{n}$ matrix, say
$(\varphi_{ij}(z))_{n\times{n}}$.

Let $F$ be a non-empty family of functions in $H^(\Omega)$. Recall
the definition of the greatest common quasi-inner divisor $\theta$
of $F$ which means that $\theta$ divides every element in $F$, and
$\theta$ is a multiple of any other common quasi-inner divisors of
$F$.
 By the same way as the proof of
Proposition 2.3.4 in \cite{20}, we have the following result.

\begin{prop}\label{24}
Let $F$ be a non-empty family of functions in $H^(\Omega)$. Then
the greatest common quasi-inner divisor $\theta$ of $F$, denoted
by $\bigwedge{F}$, exists.

\end{prop}

\begin{lem}\label{68}
 If
$\varphi=[\varphi_{ij}]_{i,j=1}^{n}$ is a function in
${H^{\infty}(\Omega,L(\C^{n}))}$, then
$D(\varphi)=\bigwedge{\varphi_{ij}}$.
\end{lem}
\begin{proof}
Let $\bigwedge{\varphi_{ij}}=\theta_{1}$. If
$\theta\in{H^{\infty}(\Omega)}$ is a quasi-inner function such
that $\theta|\varphi$, then $\varphi=\theta\phi$ for some
$\phi\in{H^{\infty}(\Omega,L(\C^{n}))}$. Thus
\begin{center}$(\varphi_{ij})_{n\times{n}}=(\theta\phi_{ij})_{n\times{n}}$\end{center} for
each $i, j=1,2,\cdot\cdot\cdot,n$. Thus
$\varphi_{ij}=\theta\phi_{ij}$ for each $i,
j=1,2,\cdot\cdot\cdot,n$. Since
$\phi\in{H^{\infty}(\Omega,L(\C^{n}))}$,
$\phi_{ij}\in{H^{\infty}(\Omega)}$ and so $\theta|{\varphi_{ij}}$
 for each $i,
j=1,2,\cdot\cdot\cdot,n$. By definition of greatest common
divisor, we can conclude that $\theta|\theta_{1}$. Clearly
$\theta_{1}|\varphi$. Thus $D(\varphi)=\theta_{1}.$
\end{proof}

We have another result similar to Theorem \ref{34}.

\begin{thm}\label{65}
For any quasi-inner functions $\theta\in{H^{\infty}(\Omega)}$ and
$\varphi\in{H^{\infty}(\Omega,L(K))}$, the following assertions
are equivalent:

(a) $\varphi{|}{\theta}$.

(b)
$\theta{H}^{\infty}(\Omega,K)\subset\varphi{H}^{\infty}(\Omega,K)$.

(c) $\theta{H}^{2}(\Omega,K)\subset\varphi{H}^{2}(\Omega,K)$.

(d) There is a $\lambda>{0}$ such that
$|\theta(z)|^{2}I_{K}\leq{\lambda^{2}\varphi(z)\varphi(z)^{\ast}}$
for any $z\in{\Omega}$.

\end{thm}
\begin{proof}
This theorem is proven by the same way as Theorem \ref{34}.
\end{proof}

\begin{prop}\label{35}
For any quasi-inner functions $\varphi$ and $\varphi^\prime$ in
$H^{\infty}(\Omega,L({\C^n}))$ and $H^{\infty}(\Omega,L(C^{n}))$
respectively, the following assertions are equivalent:

(a) $\varphi{|}_{l}{\varphi^\prime}$.

(b)
$\varphi^{\prime}{H}^{\infty}(\Omega,{\C}^{n})\subset\varphi{H}^{\infty}(\Omega,{\C^n})$.

(c)
$\varphi^{\prime}{H}^{2}(\Omega,{\C}^{n})\subset\varphi{H}^{2}(\Omega,{\C^n})$.
\end{prop}
\begin{proof}
In the same way as the proofs of Theorem \ref{34}, it is proven
that (a) and (b) are equivalent and so are (a) and (c).
\end{proof}
\begin{cor}\label{36}
Let $\varphi$ and $\varphi^\prime$ be quasi-inner functions in
$H^{\infty}(\Omega,L({\C^n}))$. Then
$\varphi{|}_{l}{\varphi^\prime}$ and
$\varphi^{\prime}{|}_{l}{\varphi}$ if and only if there exist a
function $\psi$ in

\noindent${H^{\infty}(\Omega,L(K_{1},K_{2}))}$ such that $\psi(z)$
is invertible for any $z\in{\Omega}$ with $\sup
\norm{\psi(z)^{-1}}<\infty$ and
$\varphi^{\prime}(z)=\varphi(z)\psi(z)$ for any $z\in\Omega$.
\end{cor}
\begin{proof}
By Proposition \ref{35}, we can conclude that
$\varphi{|}_{l}{\varphi^\prime}$ and
$\varphi^{\prime}{|}_{l}{\varphi}$ if and only if
$\varphi^{\prime}{H}^{2}(\Omega,{\C}^{n})=\varphi{H}^{2}(\Omega,{\C^n})$.
Then by Proposition \ref{23}, this Corollary is proven.
\end{proof}

For quasi-inner functions $\varphi$ and $\varphi^\prime$ in
${H^{\infty}(\Omega,L({\C^n}))}$, $\varphi(z)$ is a bounded
operator on the Hilbert space $\C^n$. When we consider some
relationships between $\varphi|_{l}{\varphi^\prime}$
 ($\varphi|_{r}{\varphi^\prime}$) and
$\varphi(z)|_{l}{\varphi^\prime}(z)
(\varphi(z)|_{r}{\varphi^\prime}(z)$, respectively$)$ for each
$z\in{\Omega}$, first of all, we obtain easily if
$\varphi|_{l}\varphi^{\prime}(\varphi|_{r}{\varphi^\prime}$) then
\[\varphi(z)|_{l}{\varphi^\prime}(z)(\varphi(z)|_{r}{\varphi^\prime}(z),
\texttt{ respectively})\] for each $z\in{\Omega}$.

Suppose that $z\in{\Omega}$ if and only if
$\overline{z}\in{\Omega}$. Then for
$\varphi\in{H^{\infty}(\Omega,L({\C^n}))}$, we define a new
function $\varphi^{\sim}\in{H^{\infty}(\Omega,L({\C^n}))}$ by
\begin{center}$\varphi^{\sim}(z)=\varphi(\overline{z})^{*}$\end{center}
for $z\in{\Omega}$.
\begin{lem}\label{61}
Let $\varphi_{1}$ and $\varphi_{2}$ be functions in
${H^{\infty}(\Omega,L({\C^n}))}$. Suppose that $z\in{\Omega}$ if
and only if $\overline{z}\in{\Omega}$. Then following statements
are equivalent:

(a) $\varphi_{1}|_{l}\varphi_{2}$.

(b) $\varphi_{1}^{\sim}|_{r}\varphi_{2}^{\sim}$.
\end{lem}
\begin{proof}
($\Rightarrow$) Since $\varphi_{1}|_{l}\varphi_{2}$, there is a
function $\varphi_{3}\in{H^{\infty}(\Omega,L({\C^n}))}$ such that
\begin{center}$\varphi_{2}(z)=\varphi_{1}(z)\varphi_{3}(z)$\end{center} for $z\in{\Omega}$.
By our assumption,
$\varphi_{2}(\overline{z})=\varphi_{1}(\overline{z})\varphi_{3}(\overline{z})$
for $z\in{\Omega}$. Thus
\begin{center}$\varphi_{2}(\overline{z})^{*}=\varphi_{3}(\overline{z})^{*}\varphi_{1}(\overline{z})^{*}$\end{center}
for $z\in{\Omega}$ and so
$\varphi_{1}^{\sim}|_{r}\varphi_{2}^{\sim}$.

($\Leftarrow$) By the same way as the proof of ($\Rightarrow$),
this implication can be proven.
\end{proof}
From Lemma \ref{61}, we have the following result.
\begin{cor}\label{62}
Let $\varphi_{1}$ and $\varphi_{2}$ be functions in
${H^{\infty}(\Omega,L({\C^n}))}$. Suppose that $z\in{\Omega}$ if
and only if $\overline{z}\in{\Omega}$. Then following statements
are equivalent:\vskip0.2cm

(a) $\varphi_{1}|\varphi_{2}$.

(b) $\varphi_{1}|_{r}\varphi_{2}$ and
$\varphi_{1}^{\sim}|_{r}\varphi_{2}^{\sim}$.

\end{cor}

------------------------------------------------------------------------

\bibliographystyle{amsplain}
\bibliography{xbib}
\end{document}